\documentclass[11pt]{amsart}

\usepackage{amsfonts,amssymb,stmaryrd,amscd,amsmath,latexsym,amsbsy}

\usepackage{amssymb}
\usepackage{amsfonts}
\usepackage{latexsym}

\newtheorem{theorem}{Theorem}[section]
\newtheorem{lemma}[theorem]{Lemma}
\newtheorem{proposition}[theorem]{Proposition}

\theoremstyle{definition}

\newtheorem{remark}[theorem]{Remark}

\newcommand{\Ext}{\text{Ext}}
\newcommand{\Res}{\text{Res}}

\newcommand{\ben}{\begin{enumerate}}
\newcommand{\een}{\end{enumerate}}

\newcommand{\Vect}{{\text{Vect}}}

\newcommand{\CC}{{\mathbb{C}}}
\newcommand{\ZZ}{{\mathbb{Z}}}

\newcommand{\tC}{{\widetilde{\mathcal C}}}

\theoremstyle{plain}

\newtheorem{thm}{Theorem}[section]

\newtheorem{prop}{Proposition}[section]

\newtheorem{lem}{Lemma}[section]

\newtheorem{cor}{Corollary}[section]

\newtheorem*{sol}{Solution}

\theoremstyle{definition}

\theoremstyle{remark}

\newcommand{\solu}[1]{\begin{sol}{\bf (\ref{#1})}}

\newcommand{\oA}{{\overline{A}}}

\begin{document}

\title{New deformations of group algebras of Coxeter groups, II}

\author{Pavel Etingof}
\email{etingof@math.mit.edu}

\author{Eric Rains}
\email{rains@math.ucdavis.edu}
\maketitle

\centerline{Dedicated to the 60-th birthday of Joseph Bernstein}

\section{Introduction}

This paper is a sequel of \cite{ER}. Specifically, 
let $W$ be a Coxeter group, generated by $s_i,i\in I$. 
Then, following \cite{ER}, one can define a new deformation $A_+=A_+(W)$ 
of the group algebra $\mathbb Z[W_+]$ of the group $W_+$ of even
elements in $W$. This deformation is an algebra over the ring
$R=\mathbb Z[t_{ijk}^{\pm 1}]=
\mathbb Z[\mathbb T]$ of regular functions on a certain
torus $\mathbb T$ of deformation parameters. The main result of 
\cite{ER} implies that this deformation is flat (i.e., $A_+$ is a flat
$R$-module) if and only if for every triple of indices 
$\Delta=\lbrace{i,j,k\rbrace}\subset I$ the corresponding rank 3 parabolic 
subgroup $W_\Delta\subset W$ is infinite.\footnote{To be more
precise, in \cite{ER} we work over $\mathbb C$, but the results 
routinely extend to the case of ground ring $\mathbb Z$.} 

Unfortunately, this result is not entirely satisfactory, since 
the most interesting Coxeter groups (finite, affine, hyperbolic)
typically do not satisfy the condition that all the rank 3
parabolic subgroups are infinite. For this reason, it is interesting to
generalize the theory of \cite{ER} to the case of general Coxeter
groups. Such a generalization is the main goal of this paper.

More precisely, recall that in \cite{ER} we constructed 
elements $T_{w(x)}$, labeled by $x\in W_+$, which span $A_+$
over $R$. Let $J\subset R$ be the intersection of all ideals 
$I$ in $R$ such that $T_{w(x)}$ are a basis of $A_+/IA_+$ over
$R/I$, and let $\mathbb T_f={\rm Spec}R/J$. If $k$ is a field,
then $\mathbb T_f(k)$ is the subset of points
$u\in \mathbb T(k)$ where the fiber $A_{u+}$ 
of the algebra $A_+$ is ``a flat deformation of $\mathbb Z[W_+]$''
(in the sense that the spanning elements $T_{w(x)}\in A_{u+}$ are actually a
basis). Thus we call $\mathbb T_f$ the 
{\bf flatness locus}. The main result of the paper
is an explicit characterization of the flatness locus. 
Namely, for every triple of indices $\Delta=\lbrace{i,j,k\rbrace}
\subset I$ denote by $\mathbb T_f^\Delta$ the
flatness locus for the ``rank 3 subalgebra'' 
$A_+^\Delta$ of $A_+$ generated by 
$a_{ij},a_{jk},a_{ki}$. Note that if $W_\Delta$ is
infinite, then $\mathbb T^\Delta_f=\mathbb T$. 

It is clear that $\mathbb T_f\subset
\cap_{\Delta\subset I}\mathbb T_f^\Delta$. One of our main results
is the following theorem. 

\begin{theorem}\label{mai}
One has $\mathbb T_f=\cap_{\Delta\subset I}\mathbb T_f^\Delta$.
\end{theorem}

In other words, the algebra $A_+$ is flat at $u$ if and only
if so is $A_+^\Delta$ for every $\Delta$. Obviously, this is a
generalization of the main result of \cite{ER}.

Theorem \ref{mai} reduces the problem of explicit description 
of the flatness locus $\mathbb T_f$ for any Coxeter group $W$ 
to the same problem for finite Coxeter groups of rank 3, i.e. 
$A_1\times I_{2m}$, $A_3$, $B_3$, and $H_3$. 
The even subgroups of these groups have a presentation 
by three generators $x,y,z$ of specified orders, with the
additional relation that the product of the generators is $1$. 
Therefore, the problem of finding $\mathbb T_f$ for such groups
reduces to the so called multiplicative Deligne-Simpson problem, 
which is the problem of existence of a collection of matrices 
with given spectra whose product is the identity matrix.
Thus, using known results about the Deligne-Simpson problem, we 
obtain an explicit description of $\mathbb T_f$ in the rank 3 case. 

Our description of $\mathbb T_f$ 
implies that $\mathbb T_f$ is a union of affine
subtori of $\mathbb T$ (possibly of different dimensions).  These tori are not
disjoint, however, and are quite large in number, since they
break the symmetry of the action of a large product of symmetric
groups on $\mathbb T_f$.  This suggests looking at the quotient of $\mathbb T_f$
by this permutation action, which has a nicer geometric
structure. Namely, over $\CC$, this quotient is the total space of a
vector bundle (with dimension depending on the component) over a
commutative affine algebraic group, the zero section of which parametrizes
points $u$ such that the algebra $A_{u+}$ is a twisted group algebra
of $W_+$.  Over $\ZZ$, this structure nearly carries over; it is still a
scheme over the locus of twisted group algebras cut out by linear equations
(so the fibers are affine spaces), but can fail to be locally trivial near a
nonreduced point of the base.  (Note that the base is nonreduced in
characteristic p if and only if p divides the number of components of
$H^*(W_+,\mathbb C^*)$.)

In particular, the point $1\in \mathbb T$ (at which all the
generators $a_{ij}$ of $A_+$ are unipotent) belongs to $\mathbb
T_f$; the 2-cocycle $\psi$ corresponding to 
this point is the 2-cocycle $\psi_{\rm spin}$ of the spinor
representation of $W_+$. 

Note that a special case of Theorem \ref{mai} and the explicit description of
2-cocycles in the rank 3 case yield an explicit description of
$H^2(W_+,\mathbb C^*)$. This description (which is 
an important ingredient in the proof of Theorem \ref{mai}) 
was actually known before, and is contained in the paper \cite{Bu}.

At the end of the paper, we show that the classical Iwahori-Hecke algebra is a
special case of our deformation. 
We also consider the additive versions ${\mathcal A}_{u+}$ 
of the algebras $A_{u+}$, and study their properties. 
In particular, we show that the Hilbert series of ${\mathcal A}_{0+}$ is 
$h(z)/(1+z)$, where $h(z)$ is the growth series of the Coxeter
group $W$. 

{\bf Acknowledgments.} 
It is our pleasure to dedicate this paper to Joseph Bernstein.
His work, as well as style of doing and explaining mathematics,
have been an inspiration for generations of mathematicians. 
 
The work of P.E. was partially supported by the NSF grant 
DMS-0504847 and the CRDF grant RM1-2545-MO-03. 
E.R. was supported in part by NSF grant DMS-0401387.
The authors would also like to acknowledge that at many stages of this
work they used the MAGMA package for algebraic computations
\cite{Ma}. 

\section{Preliminaries}

\subsection{Coxeter groups}
Let $I$ be a finite set, and let $\mathbb Z_{\ge 2}$ denote the set of
integers which are $\ge 2$. A Coxeter matrix
over $I$ is a collection $M$ of elements
$m_{ij}\in \mathbb Z_{\ge 2}\cup \lbrace{\infty\rbrace}$,
$i,j\in I$, $i\ne j$, such that $m_{ij}=m_{ji}$.
The rank $r$ of $M$ is, by definition, the cardinality of the set
$I$.

Let $M$ be a Coxeter matrix. Then one defines the Coxeter group
$W(M)$ by generators $s_i$, $i\in I$, and defining relations
$$
s_i^2=1,\ (s_is_j)^{m_{ij}}=1 \text{ if }m_{ij}< \infty.
$$
For brevity, from now on we will not write the dependence on $M$
explicitly, assuming that $M$ has been fixed (unless confusion is possible). 

The group $W$ has a sign character
$\xi: W\to \lbrace{\pm 1\rbrace}$ given
by $\xi(s_i)=-1$. Denote by $W_+$ the kernel of $\xi$,
i.e. the subgroup of even elements. 

Let $J\subset I$ be a subset. Let $W_J$ be the subgroup of $W$ 
generated by $s_i,i\in J$. Then $W_J$ is a Coxeter group, which
is called the parabolic subgroup of $W$ corresponding to $J$.  

Let $V$ be the reflection representation of $W$, i.e. 
the space $\mathbb C^r$ with the action of $W$ given by 
$s_i(e_j)=e_j+2{\rm cos}(\pi/m_{ij})e_i$, 
$s_i(e_i)=-e_i$. This action $\rho: W\to GL(V)$ preserves the 
(possibly degenerate) inner product on $V$ given by 
$(e_i,e_i)=1, (e_i,e_j)=-{\rm cos}(\pi/m_{ij})$. (here by
definition if $m_{ij}=\infty$ then $\pi/m_{ij}=0$).
Thus $\rho(W)\subset O(V)$, $\rho(W_+)\subset SO(V)$.

For any $x\in W(M)$, let $l(x)$ be the length of $x$. 
Let $b_n$ be the number of elements of $W(M)$ with $l(x)=n$. 
The series $h_M(z):=\sum_{n\ge 0}b_nz^n$ is called the 
growth series of $W$. It is known to be a rational function 
and can be computed explicitly (\cite{B}). 

\subsection{The algebras $A,A_+, A_u,A_{u+}$.}
Recall the definition of the algebras $A(M)$, $A_+(M)$ from \cite{ER}.

Let $t_{ijk}$ be commuting variables defined for $i\ne j$ such that
$m_{ij}<\infty$, and $k\in \mathbb Z_{m_{ij}}$
(an integer modulo $m_{ij}$), so that
$t_{ijk}=t_{ji,-k}^{-1}$. Let $R=\mathbb Z[t_{ijk}]$ be the ring of
polynomials of these variables.  (In \cite{ER}, we took coefficients in
$\mathbb C$, but for all of the results cited in this section, the same proof
works over $\mathbb Z$.)  Define the algebra $A$ generated over $R$ by
generators $s_i$ with defining relations
$$
s_i^2=1,\
\prod_{k=1}^{m_{ij}}(s_is_j-t_{ijk})=0\text{ if }m_{ij}<\infty,
$$
$$
s_pt_{ijk}=t_{jik}s_p.
$$
Let $A_+$ be the subalgebra of $A$ generated over $R$ by 
$a_{ij}:=s_is_j$. The defining relations for $A_+$
are 
$$
\prod_{k=1}^{m_{ij}}(a_{ij}-t_{ijk})=0\text{ if }m_{ij}<\infty,
$$
$$
a_{ij}a_{ji}=1, a_{ij}a_{jp}a_{pi}=1.
$$
Note that $R$ is central in $A_+$ (but not in $A$). 

Define the algebraic torus $\mathbb T={\rm Spec}R$. The group 
$\mathbb Z_2$ acts on 
$\mathbb T$ via $t_{ijk}\to t_{jik}$. If $u\in \mathbb T$ is a
point over some field, then we denote 
by $A_u,A_{u+}$ the fibers of $A,A_+$ at $u$. 
That is, $A_u=A/J_uA$, $A_{u+}=A_+/J_uA_+$, where $J_u\subset R$
is the maximal ideal corresponding to $u$. 
Note that $A_{u+}$ is an algebra, while $A_u$ is an algebra 
for any $u$ that is fixed under $\mathbb Z_2$.
In particular, we can define the ``unipotent'' algebras 
$A_1,A_{1+}$, where $1\in \mathbb T$ 
is the unit, i.e., the point where $t_{ijk}=1$ 
for all $i,j,k$. Explicitly, the algebra $A_1$ is generated by $s_i$
with defining relations 
$$
s_i^2=1,\ (s_is_j-1)^{m_{ij}}=0\ (m_{ij}<\infty)
$$
and $A_{1+}$ is generated by
$a_{ij}=s_is_j$, with defining relations
$$
a_{ij}a_{ji}=1, a_{ij}a_{jk}a_{ki}=1, (a_{ij}-1)^{m_{ij}}=0
(m_{ij}<\infty).
$$
 
Consider the subgroup $Z\subset \mathbb G_m^{r(r-1)/2}$ of collections
\[
\lbrace{z_{ij}\in \mathbb G_m,i\ne j|
z_{ij}z_{ji}=1,z_{ij}z_{jk}z_{ki}=1\rbrace}
\cong
\mathbb G_m^{r-1}.
\]
It acts on the algebra $A_+$ by rescalings: 
$a_{ij}\to z_{ij}a_{ij}$, $t_{ijk}\to z_{ij}t_{ijk}$. 
Thus the algebras $A_{u+}$ are, in effect, parametrized by 
the quotient group $\mathbb T/Z$ (in the sense that $A_{u+}$
and $A_{zu+}$ are naturally isomorphic for $z\in Z$).  

Let $x\to w(x)$ be a function assigning to every element $x\in W$ a reduced
word representing $x$. Also, for any word $w$ in letters $s_i$, let $T_w$
be the corresponding element of $A$ (or $A_u$). Recall from \cite{ER} that
$T_{w(x)},x\in W$, span $A$ as a left $R$-module, and hence span $A_u$ for
each $u$ (and similarly, $T_{w(x)}$, $x\in W_+$, span $A_+$ and $A_{u+}$).
Let $J\subset R$ be the intersection of all ideals 
$I$ in $R$ such that $T_{w(x)}$ are a basis of $A_+/IA_+$ over
$R/I$, and let $\mathbb T_f={\rm Spec}R/J$. (It is easy
to see (cf. \cite{ER}) that this property of $I$ 
is independent of the function $w(x)$). 
Thus $\mathbb T_f\subset \mathbb T$ is a closed
subscheme of $\mathbb T$, which we call the {\bf flatness locus}. 
Note that for any field $k$, $\mathbb T_f(k)$ is the set of all 
$u\in \mathbb T(k)$ such that the elements $T_{w(x)}, x\in W$ 
are a basis of $A_u$ (equivalently $T_{w(x)}$, $x\in W_+$, 
are a basis of $A_{u+}$).

Since $A_{u+}$ is unchanged if we reorder the $t_{ijk}$ for each $i<j$ (and
reorder $t_{jik}$ accordingly), we find that $\mathbb T_f$ is preserved by the
natural action of the group $\prod_{i<j} S_{m_{ij}}$ on $\mathbb T$
(the product is taken over those $i,j$ for which $m_{ij}<\infty$).  The
quotient $\mathbb T/\prod_{i<j} S_{m_{ij}}$ is the spectrum of the ring
\[
\tilde{R}:=\ZZ[e^{(k)}_{ij}\ (1\le k\le m_{ij}),1/e^{(m_{ij})}_{ij}],
\]
where
\[
e^{(k)}_{ji} = e^{(m_{ij}-k)}_{ij}/e^{(m_{ij})}_{ij},\ e^{(m_{ij})}_{ji} = 1/e^{(m_{ij})}_{ij}.
\]
This gives rise to an algebra $\tilde{A}$ with presentation
\[
s_i^2 = 1,\ (s_is_j)^{m_{ij}} + \sum_{k=1}^{m_{ij}} (-1)^k e^{(k)}_{ij}
(s_is_j)^{m_{ij}-k} = 0\ \text{ if }\ m_{ij}<\infty,\ s_p e^{(k)}_{ij} =
e^{(k)}_{ji} s_p,
\]
and a subalgebra $\tilde{A}_+$ with presentation
\[
a_{ij}^{m_{ij}} + \sum_{k=1}^{m_{ij}} (-1)^k e^{(k)}_{ij}
a_{ij}^{m_{ij}-k} = 0\ if\ m_{ij}<\infty,
\ a_{ij}a_{ji}=1,\ a_{ij}a_{jp}a_{pi}=1.
\]

Denote by $\tilde A_v$, $\tilde A_{v+}$ the fiber of 
$\tilde A,\tilde A_+$, respectively, at the point $v\in 
\mathbb T/\prod_{i<j} S_{m_{ij}}$. Let $\pi: 
\mathbb T\to \mathbb T/\prod_{i<j} S_{m_{ij}}$ be the natural
projection. Note that $A_{u+}\cong \tilde{A}_{\pi(u)+}$ and $A_u\cong
\tilde{A}_{\pi(u)}$, and thus $A_u$ has a natural algebra
structure not only when $u$ is $\mathbb Z_2$-invariant, 
but also in the more general case when only $\pi(u)$ is 
$\mathbb Z_2$-invariant.  In any event, we observe that the flatness locus
associated to $\tilde{A}$ is 
$\tilde{\mathbb T}_f:=\pi(\mathbb T_f)$, 
since the coefficients of the
relations of $A_{u+}$ are functions on $\tilde{T}$; 
in particular, this holds
scheme-theoretically, not just point-wise.

\section{Schur multipliers} 

Let $\mathbb T_f^S$ be the set of all complex points of $u\in \mathbb T_f$ 
for which 
$$
\prod_k(z-t_{ijk})=z^{m_{ij}}+(-1)^{m_{ij}}t_{ij}
$$
for some complex numbers $t_{ij}$ ($t_{ij}=t_{ji}^{-1}$),
defined when $m_{ij}<\infty$.  
Denote the set of the corresponding 
collections of numbers $t_{ij}$ by $\Theta$
(clearly, the isomorphism class of $A_{u+}$, $u\in \mathbb T_f^S$, depends
only on the image of $u$ in $\Theta$, which we denote by $t=t_u$). 
Then the polynomial relations for $A_{u+}$ take the form
$a_{ij}^{m_{ij}}=(-1)^{m_{ij}+1}t_{ij}$. Let $[x]:=T_{w(x)}$ 
Then in the algebra $A_{u+}$ we
have the following multiplication table:
$[x][y]=\psi_t(x,y)[xy]$, where $\psi_t(x,y)$ are some scalars. 
Indeed, the scalar $\psi_t(x,y)$ is $\pm$ a product of $t_{ij}$ 
which appears when we express $w(xy)$ via $w(x)w(y)$ using the
braid relations. 
Note that the braid relations written in terms of $a_{ij}\in A_{u+}$  
 look like 
$$
a_{ij}^{m_{ij}/2}=-t_{ij}a_{ji}^{m_{ij}/2},
$$
$$ 
a_{pj}a_{ij}^{(m_{ij}-2)/2}a_{iq}=-t_{ij}a_{pi}a_{ji}^{(m_{ij}-2)/2}a_{jq}
$$ 
if $m_{ij}$ is even, and 
$$
a_{ij}^{(m_{ij}-1)/2}a_{ip}=
t_{ij}a_{ji}^{(m_{ij}-1)/2}a_{jp},
$$  
$$
a_{pi}a_{ij}^{(m_{ij}-1)/2}=
t_{ij}a_{pj}a_{ji}^{(m_{ij}-1)/2}
$$
if $m_{ij}$ is odd. 

Thus $\psi_t(x,y)$ is a Schur multiplier (a 2-cocycle 
of $W_+$ with coefficients in $\mathbb C^*$), and $A_{u+}$ 
is the twisted group algebra $\mathbb C_{\psi_t}[W_+]$. It is obvious that 
if the function $w(x)$ is changed, then the 2-cocycle $\psi_t$
changes by a coboundary. Thus we have a canonical map 
$\psi: \Theta\to H^2(W_+,\mathbb C^*)$ defined by $t\to \mathbb
\psi_t$. 

\begin{proposition}\label{prop:Schur_flatness}
A system of numbers $t_{ij}$ is in $\Theta$ if and only if
whenever $W_{ijk}$
is
finite,
\[
t_{ij}^{[W_{ijk}:W_{ij}]}
t_{jk}^{[W_{ijk}:W_{jk}]}
t_{ki}^{[W_{ijk}:W_{ki}]}=1.
\]
\end{proposition}

\begin{proof} By definition, $\{t_{ij}\}\in\Theta$ iff the
associated algebra $A_{u^+}$ is flat.  Since this has the form of
a twisted
group algebra, this may be restated group-theoretically as
follows.  Let
$F$ be the free group on $r$ elements $s_1$,\dots,$s_r$, and let
$G:=\CC^*\ltimes F$, where each generator of $F$ acts on $\CC^*$ by
$x\mapsto
x^{-1}$.  Finally, let $N$ be the normal subgroup of $G$
generated by the
elements
\begin{align}
s_i^2,&&1\le i\le r\\
(-1)^{m_{ij}+1}t_{ij}(s_is_j)^{-m_{ij}},&&1\le i<j\le r.
\end{align}
Then $t_{ij}\in \Theta$ iff $N$ is disjoint from $\CC^*$, or
equivalently
iff there exists a (noncentral) extension of $W$ by $\CC^*$ having
preimages
of the generators satisfying the above relations.  The main
theorem
of \cite{Bu} gives necessary and sufficient conditions for
extensions of
Coxeter groups by arbitrary $\ZZ$-modules; specializing
Burichenko's
condition to the present case gives the desired result.
\end{proof}

Note that the group $Z$ acts on $\Theta$ via $t_{ij}\to
t_{ij}z_{ij}^{m_{ij}}$.  It is easy to see that the map $\psi$ descends to
a map $\psi: \Theta/Z\to H^2(W_+,\mathbb C^*)$, which is injective.  Let us
show that this map is in fact a bijection.  To do so, we will construct the
inverse map $\eta$.  Let $\psi$ be a Schur multiplier for $W_+$, and $\mathbb
C_\psi[W_+]$ be the corresponding twisted group algebra.  Let $a_{ij}^{\rm
  gr}$ be the elements $a_{ij}$ in the group $W_+$ as opposed to the
algebra $\mathbb C_\psi[W_+]$.  Let $0\in I$ be an element. Define
$a_{0j}:=[a_{0j}^{\rm gr}]$ for $j\ne 0$, $a_{j0}:=a_{0j}^{-1}$, and
$a_{ij}:=a_{0i}^{-1}a_{0j}$ for $i\ne 0,j\ne 0$.  Then
$a_{ij}^{m_{ij}}=t_{ij}$, where $t_{ij}=t_{ij}(\psi)$.  Thus we have
attached to $\psi$ an element $t$ of $\Theta$ (and thus of $\Theta/Z$); we
set $t=\eta(\psi)$. It is easy to see that $\eta$ is the desired inverse.

As a by-product, we see that $\Theta/Z$ is a group under
multiplication, and $\psi,\eta$ are group isomorphisms. 

In the sequel, a particular class in $H^2(W_+,\mathbb C^*)$
will be especially important. It is the 
Schur multiplier $\psi_{\rm Spin}$ afforded by 
the pullback of a spinor representation
of (the reductive part of) $SO(V)$ under the homomorphism  
$\rho: W_+\to SO(V)$. 

\begin{lemma}\label{spin} 
If $\psi=\psi_{\rm spin}$ then $\eta(\psi)$ 
is the class of $t=1$, (i.e. $t_{ij}=1$ for all $i,j$). 
\end{lemma}

\begin{proof} The representation 
$\rho: W\to O(V)$ lifts to a projective representation 
$\hat\rho: W\to {\rm Pin}(V)$. This representation is defined by the
formula $\hat\rho(s_i)=e_i$. We have $e_i^2=1$, and
$e_ie_j+e_je_i=
-2{\rm cos}(\pi/m_{ij})$. For $m_{ij}<\infty$, 
let $Y_{ij}$ be the 2-dimensional spin
representation of the Clifford algebra generated by $e_i,e_j$. 
Then it follows from the last equation that the trace of $e_ie_j$ in this
representation is $-2{\rm cos}(\pi/m_{ij})$. Thus the eigenvalues
of  $e_ie_j$ in $Y_{ij}$ are $-e^{\pm\pi\sqrt{-1}/m_{ij}}$, and
hence $(e_ie_j)^{m_{ij}}=(-1)^{m_{ij}+1}$, as desired.
\end{proof}

\section{The flatness locus for finite triangles}

In \cite{ER}, it was shown that if all triangular (i.e., rank 3 parabolic)
subgroups of $W$ are infinite, then $\mathbb{T}_f=\mathbb{T}$; that is, all
parameters are flat.  This suggests that in general, finite triangles will
play a particularly important role.

It is a classical result that the triangle group with exponents $p$, $q$,
$r$ is finite if and only if $1/p+1/q+1/r>1$, and thus the only finite
triangle groups are the infinite family
\[
\langle x,y,z|x^2=y^2=z^n=xyz=1\rangle,\quad 2\le n
\]
and the three sporadic examples
\[
\langle x,y,z|x^2=y^3=z^n=xyz=1\rangle,\quad 3\le n\le 5.
\]
We may thus proceed by a case-by-case analysis.

Since we are dealing with a finite algebra, the flatness condition at a
point $u$ is equivalent to the requirement that $\dim(A_{u+})=|W^+|$.  If
$A_{u+}$ is semisimple, this can be verified by exhibiting sufficiently
many irreducible representations (with dimensions satisfying $\sum d^2 =
|W^+|$).  This suffices in particular to determine any component of
$\mathbb{T}_f$ for which the generic point is semisimple, or equivalently
any component that contains {\em any} semisimple point.  In particular, the
group algebra itself is such a point, and thus gives rise to a large
portion of the flatness locus.  The representation theory of the group
algebra is obviously relevant, and we must therefore distinguish in the
infinite family between the cases $n$ even and odd.

\begin{lem}\label{lem:flat_triangle_group}
Suppose $n\ge 2$ is even, and consider the algebra $A_{u+}$ with
$m_{12}=m_{23}=2$, $m_{13}=n$.  If
\begin{align}
t_{122}t_{232}&=t_{13n},\notag\\
t_{121}t_{231}&=t_{13n},\notag\\
t_{122}t_{231}&=t_{13(n/2)},\notag\\
t_{121}t_{232}&=t_{13(n/2)},\notag\\
t_{121}t_{122}t_{231}t_{232}&=t_{13i}t_{13(n-i)},\quad 1\le i<n/2\notag
\end{align}
then $A_{u+}$ is flat.  Similarly, if $n\ge 3$ is odd, $A_{u+}$ is flat if
\begin{align}
t_{122}t_{232}&=t_{13n}\notag\\
t_{121}t_{231}&=t_{13n}\notag\\
t_{121}t_{122}t_{231}t_{232}&=t_{13i}t_{13(n-i)},\quad 1\le i<n/2\notag.
\end{align}
For the exceptional cases $m_{12}=2$, $m_{23}=3$, $m_{13}=n$, the algebra
is flat if $n=3$ and
\begin{align}
t_{122}t_{233}&=t_{133}\notag\\
t_{122}t_{231}&=t_{131}\notag\\
t_{122}t_{232}&=t_{132}\notag\\
t_{121}^2t_{122}t_{231}t_{232}t_{233} &= t_{131}t_{132}t_{133};\notag
\end{align}
if $n=4$ and 
\begin{align}
t_{122}t_{233}&=t_{134}\notag\\
t_{121}t_{233}&=t_{132}\notag\\
t_{121}t_{122}t_{231}t_{232}&=t_{132}t_{134}\notag\\
t_{121}t_{122}^2 t_{231}t_{232}t_{233} &= t_{131}t_{132}t_{133}\notag\\
t_{121}^2t_{122} t_{231}t_{232}t_{233} &= t_{131}t_{133}t_{134};\notag
\end{align}
or if $n=5$ and
\begin{align}
t_{122}t_{233}&=t_{135}\notag\\
t_{121}t_{122}^2 t_{231}t_{232}t_{233} &= t_{131}t_{134}t_{135}\notag\\
t_{121}t_{122}^2 t_{231}t_{232}t_{233} &= t_{132}t_{133}t_{135}\notag\\
t_{121}^2t_{122}^2 t_{231}t_{232}t_{233}^2 &=
t_{131}t_{132}t_{133}t_{134}\notag\\
t_{121}^2 t_{122}^3 t_{231}^2t_{232}^2t_{233} &=
t_{131}t_{132}t_{133}t_{134}t_{135}.\notag
\end{align}
\end{lem}

\begin{proof}
Note that in each case, the equations cut out a subgroup scheme of $\mathbb
T$, and each of the (at most 2) components of the subgroup is smooth and
has points in characteristic 0.  Since $\mathbb T_f$ is a closed
subscheme of $\mathbb T$, we may thus work over $\mathbb C$.

It then suffices to show that the irreducible representations of $W_+$
deform, for generic solutions of the above equations; more precisely, we
claim that each equation gives the condition for such deformation.  For
instance, to deform the $5$-dimensional representation of the $(2,3,5)$
group, we must have matrices $A$ and $B$ such that $A$ has eigenvalues
$t_{121}$ and $t_{122}$ with multiplicities 2 and 3, $B$ has eigenvalues
$t_{231}$, $t_{232}$, $t_{233}$ with multiplicities $2$, $2$, and $1$, and
$AB$ has eigenvalues $t_{13i}$ with multiplicity 1.  This is an example of
a rigid multiplicative Deligne-Simpson problem, and it follows from
\cite{S} that such a solution exists (and is unique up to conjugacy) if and
only if the determinants multiply up correctly.
\end{proof}

\begin{remark}
The rigid Deligne-Simpson problems that arise above and in Lemma
\ref{lem:flat_triangle_spin} are those with generically diagonalizable
matrices with eigenvalue multiplicities from the list
\begin{align}
&(1,1,1), (11,11,11), (21,111,111),\notag\\
&(22,211,1111), (32,221,11111), (33,222,21111).\notag
\end{align}
\end{remark}

Now, in most of these cases, there is another way to structure the
representations to obtain the correct dimensions; the point is that aside
from the case $(2,2,n)$, $n$ odd, the group has a nontrivial central
extension (the spin extension discussed above), giving rise to a
nontrivially twisted group algebra.

\begin{lem}\label{lem:flat_triangle_spin}
Suppose $n\ge 2$ is even, and consider the algebra $A_{u+}$ with
$m_{12}=m_{23}=2$, $m_{13}=n$.  If
\[
t_{121}t_{122}t_{231}t_{232} = t_{13i}t_{13(n+1-i)},\quad 1\le i\le n/2,
\]
then $A_{u+}$ is flat.
Similarly, the algebra is flat if $m_{12}=2$, $m_{23}=3$, $m_{13}=3$, and
\begin{align}
t_{121}t_{122}t_{231}t_{232} &= t_{131}t_{132}\notag\\
t_{121}t_{122}t_{231}t_{233} &= t_{131}t_{133}\notag\\
t_{121}t_{122}t_{232}t_{233} &= t_{132}t_{133};\notag
\end{align}
if $m_{12}=2$, $m_{23}=3$, $m_{13}=4$, and
\begin{align}
t_{121}t_{122}t_{231}t_{233}&= t_{131}t_{134}\notag\\
t_{121}t_{122}t_{231}t_{233}&= t_{132}t_{133}\notag\\
t_{121}^2t_{122}^2t_{231}t_{232}^2t_{233} &= t_{131}t_{132}t_{133}t_{134}\notag;
\end{align}
or if $m_{12}=2$, $m_{23}=3$, $m_{13}=5$, and
\begin{align}
t_{121}t_{122}t_{231}t_{232} &= t_{131}t_{135}\notag\\
t_{121}t_{122}t_{231}t_{232} &= t_{132}t_{133}\notag\\
t_{121}^2t_{122}^2t_{231}t_{232}^2t_{233} &=
t_{131}t_{132}t_{134}t_{135}\notag\\
t_{121}^3t_{122}^3t_{231}^2t_{232}^2t_{233}^2 &=
t_{131}t_{132}t_{133}^2t_{134}t_{135}.\notag
\end{align}
\end{lem}

\begin{proof}
The proof is analogous to the proof of the previous lemma. 
\end{proof}

We can of course obtain other tori inside the flatness locus by permuting
the $t_{ijk}$ for each $ij$; the claim, however, is that any flat choice of
parameters lies inside one of these tori.  We thus wish to show that $\mathbb
T_f$ is the union of a quite large number of tori; the resulting equations
would therefore be extremely complicated in general.  However, if we mod
out by the natural permutation action, we obtain a much simpler expression
for the ideal.  In other words, we need to consider the algebra
$\tilde{A}_+$, or equivalently the algebra with relations
\begin{align}
abc-1
=
a^2-\alpha_1 a+\alpha_2
=
b^2-\beta_1 b+\beta_2
&{}=0,\notag\\
c^n-\gamma_1 c^{n-1}+\gamma_2 c^{n-2}-\cdots + (-1)^n \gamma_n
&{}= 0,
\notag
\end{align}
replacing the polynomial $b^2-\beta_1 b+\beta_2$ 
with $b^3-\beta_1 b^2+\beta_2 b-\beta_3$ in
the exceptional cases.

\begin{thm}
The flatness locus for the algebra $\tilde A_+$ 
is defined by the following equations. 
For $(2,2,n)$, $n\ge 2$ even:
\begin{align}
(\alpha_2 \beta_2)^n \gamma_n^2 &= 1\notag\\
(\alpha_2\beta_2)^k \gamma_n \gamma_k &= \gamma_{n-k}\notag\\
(\alpha_2 \beta_2)^{n/2}\gamma_n \alpha_1 &= \alpha_1\notag\\
(\alpha_2 \beta_2)^{n/2}\gamma_n \beta_1 &= \beta_1;\notag
\end{align}
for $(2,2,n)$, $n\ge 3$ odd:
\begin{align}
(\alpha_2 \beta_2)^n \gamma_n^2 &= 1\notag\\
(\alpha_2\beta_2)^k \gamma_n \gamma_k &= \gamma_{n-k}; \notag\\
\alpha_2^{(n-1)/2} \beta_2^{(n+1)/2} \gamma_n \alpha_1 &= \beta_1;\notag
\end{align}
for $(2,3,3)$:
\begin{align}
\alpha_2^6 \beta_3^4 \gamma_3^4 &= 1\notag\\
\alpha_2^3\beta_3^2\gamma_3^2\alpha_1 &= \alpha_1\notag\\
\alpha_2^2 \beta_3^2 \gamma_3 \beta_1 &= \gamma_2\notag\\
\alpha_2^2 \beta_3 \gamma_3^2 \gamma_1 &= \beta_2\notag;
\end{align}
for $(2,3,4)$:
\begin{align}
\alpha_2^{12} \beta_3^8 \gamma_4^6 &= 1\notag\\
\alpha_2^6 \beta_3^4 \gamma_4^3 \alpha_1 &= \alpha_1\notag\\
\alpha_2^6 \beta_3^4 \gamma_4^3 \gamma_2 &= \gamma_2\notag\\
\alpha_2^4 \beta_3^3 \gamma_4^2 \beta_1 &= \beta_2\notag\\
\alpha_2^3 \beta_3^2 \gamma_4^2 \gamma_1 &= \gamma_3;\notag
\end{align}
and for $(2,3,5)$:
\begin{align}
\alpha_2^{30}\beta_3^{20}\gamma_5^{12} &= 1\notag\\
\alpha_2^{15}\beta_3^{10}\gamma_5^6\alpha_1 &= \alpha_1\notag\\
\alpha_2^{10}\beta_3^7\gamma_5^4 \beta_1 &= \beta_2\notag\\
\alpha_2^6\beta_3^4 \gamma_5^3 \gamma_1 &= \gamma_4\notag\\
\alpha_2^{12}\beta_3^8 \gamma_5^5 \gamma_2 &= \gamma_3\notag.
\end{align}
\end{thm}

\begin{proof}
To show that the equations are sufficient for flatness, we note
that over $\CC$ the above equations imply simple multiplicative
relations between the roots of the minimal polynomials (e.g., the
fact that if the polynomial $p$ satisfies $p(x)=(-1)^n x^n p(1/x)$,
then its roots multiply pairwise to 1), which in turn imply flatness,
using the appropriate lemma.  In particular, since flatness is a
closed condition, this proves sufficiency over $\ZZ$.

Now let us prove the necessity of the equations. 
In the exceptional cases, we can simply compute (using
MAGMA) a noncommutative Gr\"obner 
basis for $\tilde{A}_+$ over $R$, ignoring any reduced
$S$-polynomials which vanish in the group case; it turns out that the
remaining $S$-polynomials have invertible leading coefficients once
reduced.  We thus find that every element of the Gr\"obner basis in the
group case can be lifted to a monic relation in $\tilde{A}_+$, and thus
every remaining $S$-polynomial must reduce to 0 modulo those relations.
The coefficients of such reductions thus give equations for the flatness
locus.

For the infinite family, essentially the same idea holds; since in that
case, the computations must be done by hand, it is more convenient to work
with $\tilde{A}$, and stop as soon as the above equations are obtained.
\end{proof}

\begin{cor}
In each case $\mathbb T_f$ is the union of the relevant schemes from Lemmas
\ref{lem:flat_triangle_group} and \ref{lem:flat_triangle_spin}.
\end{cor}

\begin{remark}\label{bund} 
The key observation to make about the above equations is that they are
linear in the middle coefficients, and thus the fibers of the map from the
flatness locus to the locus of twisted group algebras are connected.
Moreover, in characteristic not 2, the corresponding vector spaces can be
given a basis depending only on the component; that is, the flatness locus
is a vector bundle over the locus of twisted group algebras.  Note that the
fiber over a nonreduced point will still be an affine space, but such a
point fails to have even a formal neighborhood over which the scheme
remains an affine space.
\end{remark}

Now consider the case of Coxeter groups of rank higher than 3. 
For any triple of distinct indices
$\Delta=\lbrace{i,j,k\rbrace}\subset I$, let $M_\Delta$ be the
corresponding 3 by 3 submatrix of $M$, and $W_\Delta=W(M_\Delta)$. Let 
$\mathbb T^\Delta$ be the corresponding torus of parameters;
we have a natural projection $p=p_\Delta:\mathbb T\to \mathbb T^\Delta$. 

Let $A_+^\Delta:=A_+(M_\Delta)$. 
Let $J_\Delta\subset R$ be the intersection of all ideals 
$I$ in $R$ such that $T_{w(x)},x\in W_{\Delta+}$, 
are a basis of $A^\Delta_+/IA^\Delta_+$ over
$R/I$, and let $\mathbb T_f^\Delta={\rm Spec}R/J_\Delta$
(a closed subscheme of $\mathbb T$). If $k$ is a field 
then $\mathbb T_f^\Delta(k)\subset \mathbb T(k)$ be the locus of 
points $t\in \mathbb T(k)$ such that $T_{w(x)}$, $x\in W_{\Delta +}$
are a basis of $A_{p(t)+}^\Delta$. Also, let $\tilde{\Bbb T}_f^\Delta=
{\Bbb T}_f^\Delta/\prod_{i<j}S_{m_{ij}}$. 

Since the equations on the middle coefficients are all of the form $p(x)=a
q(bx)$ or $p(x)= a x^n q(b/x)$ where $p$, $q$ are minimal polynomials of
the generators and $a$, $b$ are monomials in the constant coefficients, it
follows that similar comments to Remark \ref{bund} 
apply to the scheme $\cap_{\Delta\subset I}
\tilde{\mathbb T}^{\Delta}_f$.  Indeed, the equations on
the middle coefficients are still linear, and the only way the dimension
can depend on the constant coefficients is via equations of the form
$p(x)=p(ax)$ with $a^n=1$ or $p(x)=p(0)^{-1} x^n p(a/x)$ with $a^n=p(0)^2$.
Thus in sufficiently large characteristic (e.g., $p>m_{ij}$ whenever
$m_{ij}<\infty$), the intersection $\cap_{\Delta\subset I}\tilde{\mathbb
  T}^{\Delta}_f$ is still a vector bundle.

More generally, it follows from the above discussion that over a
sufficiently large cyclotomic ring (so the twisted group locus breaks into
smooth components), the scheme $\cap_{\Delta\subset I} \tilde{\mathbb
  T}^\Delta_f$ can be expressed as a union of smooth components, one for
each component of the twisted group locus.  Since the polynomials
associated to a point on the twisted group locus have distinct roots in a
suitable cyclotomic ring, the scheme $\cap_{\Delta\subset I} \mathbb
T^\Delta_f$ is also a union of smooth components, a total of $|\prod_{i<j}
S_{m_{ij}}|$ for each component of the twisted group locus.

\section{The structure of the flatness locus around a 2-cocycle}

\subsection{The formal PBW theorem}

Let $\mathbb T_f^*=\cap_\Delta \mathbb T_f^\Delta$.  It is clear that
$\mathbb T_f\subset \mathbb T_f^*$.  Moreover, by Proposition
\ref{prop:Schur_flatness} and the explicit computations of $\mathbb
T_f^\Delta$ in the previous section, it follows that every component of
$\mathbb T_f^*$ is integral and contains a complex $2$-cocycle $\psi\in
\Theta$.  Theorem \ref{mai} will thus follow if we manage to show that
every point of $\mathbb T_f^*$ in some neighborhood of $\psi$ is actually
in $\mathbb T_f$.

The main result of this section, which accomplishes 
this, is Theorem \ref{formal pbw} below. Note that it is enough 
for us to work over $\mathbb C$; we will do so, and 
for simplicity will use the same notation for the complex
counterparts of the $\mathbb Z$-objects considered before. 

\begin{theorem}\label{formal pbw}
Let $\psi\in \Theta$.
Then there exists a formal neighborhood $U$ of $\psi$ in $\mathbb T$ such
that $\mathbb T_f\cap U=\mathbb T_f^*\cap U$.  
\end{theorem}


The rest of the section is devoted to the proof 
of Theorem \ref{formal pbw}. 

Let $\hat R,\hat A_+$ be the completions of $R,A_+$ near the point $\psi$.
Let $J\subset \hat R$ be the completion of the ideal of $\mathbb T_f^*$. Our
job is to show that the elements $T_{w(x)}$, $x\in W_+$, form a
(topological) basis of the algebra $\hat A_+/J\hat A_+$ over $\hat R/J$.
This will be done using the theory of sheaves on posets (which is very
closely related to the constructible sheaves used in \cite{ER}, but is 
more elementary). We let $t_\psi=t$ be the collection of numbers $t_{ij}$
attached to $\psi$ as above.

\subsection{The cell complex $\Sigma$ attached to a Coxeter group}

Let us recall (see e.g. \cite{Da}) that to any Coxeter group $W=W(M)$ one
can attach a regular cell complex $\Sigma=\Sigma(M)$ as follows.  The cells
of dimension $n$ in $\Sigma$ are labelled by cosets in $W$ of a finite
parabolic subgroup $W_Q$ corresponding to an $n$-element subset $Q\subset
I$. The boundary of the cell $\sigma\in W/W_Q$ consists of the cells
$\sigma'\in W/W_{Q'}$, $Q'\subset Q$, such that $\sigma'\subset
\sigma$. Also, for every $Q$ we have a canonical cell $\sigma_Q\in W/W_Q$
containing the identity element of $W$.

\begin{theorem}\label{contra} (see e.g. \cite{Da})
The cell complex $\Sigma$ is regular and contractible. 
The group $W$ acts on $\Sigma$ properly discontinuously. 
\end{theorem}

Let $\Sigma_n$ denote the $n$-skeleton of $\Sigma$. 
It follows from Theorem \ref{contra} that 
$H^i(\Sigma_n,\mathbb C)=0$ for $i<n$. 

\subsection{Equivariant sheaves on posets}

Given a poset $P$, a {\em sheaf} on $P$ is a functor from $P$ to
$\Vect$; that is, an assignment of a (complex) vector space $S(x)$ for
$x\in P$, together with maps $f_{xy}:S(y)\to S(x)$ for $y\subset x$,
satisfying $f_{xy}\circ f_{yz}=f_{xz}$.  If the maps $f_{xy}$ are
isomorphisms, we say that $S$ is a local system.  If $V$ is a vector space,
then there is a local system (the {\em constant sheaf}) with $S(x)\equiv V$
and $f_{xy}\equiv 1$; we will denote this sheaf by $V$.

\begin{prop} (see e.g. \cite{Yu})
There is a canonical isomorphism
\[
\Ext^*_P(\CC,\CC)\cong H^*(\Delta(P)),
\]
where $\Delta(P)$ is the order complex of $P$, and $\Ext^*_P$ is taken in
the category of sheaves on $P$.
\end{prop}

In particular, if $P$ is the poset of closed cells of a regular CW complex
$\Sigma$, then $\Delta(P)$ is homeomorphic to the barycentric subdivision
of $\Sigma$, and we have $\Ext^*_P(\CC,\CC)\cong H^*(\Sigma)$.

Now, suppose the group $G$ acts on $P$ (preserving order), and let $\psi$
be a $2$-cocycle of $G$.  We can then define a {\em $\psi$-equivariant
  sheaf} to be a sheaf on $P$ together with maps
\[
\rho(g):S(x)\to S(gx)
\]
compatible with the sheaf maps and satisfying
\[
\rho(g)\rho(h) = \psi(g,h) \rho(gh).
\]
Note that for each $x\in P$, $S(x)$ is a $\CC_\psi[G_x]$-module, and the
maps $f_{xy}$ are $\CC_\psi[G_x\cap G_y]$-linear.

The forgetful functor $\pi^*$ from $\psi$-equivariant sheaves to ordinary
sheaves has a natural adjoint, which associates to the sheaf $V$ the sheaf
\[
(\pi_! V)_x = \bigoplus_{g\in G} V_{g^{-1} x}
\]
such that
\[
\rho(g):(\pi_! V)_x\to (\pi_! V)_{gx}
\]
acts on $V_{h^{-1}x}\cong V_{(g h)^{-1} gx}$ as multiplication by
$\psi(g,h)$.

\subsection{Equivariant sheaves on the 3-skeleton of $\Sigma$}

In our case, we are primarily interested in $\psi$-equivariant sheaves on
$\Sigma_3$, where $\psi$ is a 2-cocycle of $W_+$. 
To describe them explicitly, we first need to understand the structure of 
$Y_3:=\Sigma_3/W_+$. It consists of the following cells

0-cells: $N, S$ (north and south pole).

1-cells: $e_i, i\in I$, connnecting $N$ and $S$. 

2-cells: $D_{ij}$, $i<j$, a disk whose boundary is the union of $e_i$
and $e_j$, if $m_{ij}<\infty$ (for convenience we 
fix an identification $I=\lbrace{1,...,r\rbrace}$). 
The center of $D_{ij}$ has isotropy group
$W_{ij+}=\mathbb Z_{m_{ij}}$. 

3-cells: $P_{ijk}$, $i<j<k$, a 3-ball whose boundary is the sphere 
made up by $D_{ij},D_{jk},D_{ik}$, if $W_{ijk}$ is a finite group. The center 
of $P_{ijk}$ has isotropy group $W_{ijk+}$, and 
the points on the segments connecting the center of $P_{ijk}$
with the centers of $D_{ij},D_{jk},D_{ik}$ have isotropy groups
$W_{ij+}=\mathbb Z_{m_{ij}}$, 
$W_{jk+}=\mathbb Z_{m_{jk}}$, 
$W_{ik+}=\mathbb Z_{m_{ik}}$.

{\bf Remark.} More precisely, the 3-cells are of the form 
$B^3/W_{ijk+}$, where $B^3$ is the 3-ball.
But the quotient $B^3/G$ for $G\subset SO(3,\mathbb R)$ 
is always isomorphic to $B^3$ as a topological space,
since $G\subset PSL_2(\mathbb C)$ and hence $S^2/G=\mathbb CP^2/G=\mathbb
CP^1=S^2$. 
 
Thus, $\psi$-equivariant sheaves on
$\Sigma_3$ can be
specified by the following data:
\begin{itemize}
\item[1.] Spaces $V_N$, $V_S$
\item[2.] A space $V_i$ for all $i\in I$, and maps $f_{Ni}:V_N\to V_i$,
  $f_{Si}:V_S\to V_i$.
\item[3.] A $\CC_\psi[W_{ij+}]$-module $V_{ij}=V_{ji}$ for $i\ne j\in I$, and
  maps $h_{ij}:V_i\to V_{ij}$.
\item[4.] A $\CC_\psi[W_{ijk+}]$-module $V_{ijk}=V_{ikj}=\cdots$ for
  $i,j,k\in I$ distinct, and $\CC_\psi[W_{ij+}]$-linear maps
\[
f_{ij;k}:V_{ij}\to \Res^{W_{ijk+}}_{W_{ij+}} V_{ijk}
\]
\end{itemize}
subject to the obvious compatibility relations on the maps, and the
omission of spaces $V_{ij}$ and $V_{ijk}$ corresponding to infinite
parabolic subgroups.

Let ${\mathcal C}$ denote the category of sheaves on the cell poset of
$\Sigma_3$, and let ${\mathcal C}_\psi$ denote the corresponding category of
$\psi$-equivariant sheaves.  Finally, let ${\bf M}$ be the $\psi$-equivariant
local system $\pi_!(\CC)$.

\begin{lem}\label{cohom}
One has $\Ext^j_{{\mathcal C}_\psi}({\bf M},{\bf M}) = 0 $ for $j=1,2$.
\end{lem}

\begin{proof}
$$
{\rm Ext}^j_{{\mathcal C_\psi}}({\mathbf M},{\mathbf M})=
{\rm Ext}^j_{{\mathcal C_\psi}}(\pi_!\CC,\pi_!\CC)=
$$
$$
{\rm Ext}^j_{{\mathcal C}}(\CC,\pi^*\pi_!\CC)=\CC_\pi[W_+]\otimes 
{\rm Ext}^j_{{\mathcal C}}(\CC,\CC)=
$$
$$
\CC_\pi[W_+]\otimes H^j(\Sigma_3;\CC) = 0.
$$
\end{proof}

\subsection{The algebra $B$}

The category ${\mathcal C}_\psi$ is abelian, and is thus equivalent to the
category of modules over some algebra, which can in fact be read off from
the above description of the data determining a $\psi$-equivariant sheaf.
In particular, the representations of $B$ correspond to sums
\[
V_N\oplus V_S\oplus \bigoplus_i V_i\oplus \bigoplus_{i,j} V_{ij} \oplus
\bigoplus_{i,j,k}V_{ijk}
\]
as above.

Now, let $\tau=(\tau_{ijk})$, $k\in \ZZ_{m_{ij}}$, be a collection of formal
parameters.  If we replace each algebra $\CC_\psi[W_{ij+}]$,
$\CC_\psi[W_{ijk+}]$ by the corresponding deformation $A_{ij+}$,
$A_{ijk+}$, we obtain a deformed algebra $B(\tau)$, over
$\hat{R}:=\CC[[\tau]]$.  Let $J$ be the ideal of the formal neighborhood of
$\mathbb T^*_f$ in $\hat{R}$, and let $B_J(\tau)=B(\tau)/JB(\tau)$ be the
corresponding algebra over $\hat{R}_J:=\hat{R}/J$.

\begin{prop}
For any Coxeter matrix $M$, the algebra $B_J(\tau)$ is a flat deformation
of $B$ over the formal neighborhood of $t=\eta(\psi)$ in $\mathbb T^*_f$.
\end{prop}

\begin{proof}
For any set $i,j$ of distinct indices, let $p_{ij}$ be the idempotent of
$B$ which acts by $1$ on $V_N,V_S,V_i,V_j$, on $V_{ij}$ if
$m_{ij}<\infty$, and by $0$ on all the other spaces.  Also, for any
triangle $i,j,k$ such that $W_{ijk}$ is finite, let $p_{ijk}$ be the
idempotent of $B$ which acts by $1$ on $V_N,V_S,V_i,V_j,V_k$,
$V_{ij},V_{jk},V_{ik}$, $V_{ijk}$, and by $0$ on all the other spaces.

Let $N$ be the direct sum of right $B$-modules $p_{ijk}B$ and
$p_{ij}B$. Then $N$ is clearly a faithful right $B$-module, so it suffices
to show that the modules $p_{ijk}B$ and $p_{ij}B$ can be deformed to right
$B_J(\tau)$-modules.  We focus on $p_{ijk}B$; the case of 
$p_{ij}B$ is analogous (and
was in any case treated in \cite{ER}, Proposition 3.10).

Clearly, $B$ preserves the kernel of $p_{ijk}$, so $p_{ijk}B=p_{ijk}B
p_{ijk}$, and thus $B_{ijk}:=p_{ijk}B$ is a unital algebra with unit
$p_{ijk}$.  Moreover, both $p_{ijk}B$ and $p_{ijk}B(\tau)$ are simply the
analogues of $B$ for the rank 3 subsystem $\{i,j,k\}$.  We may thus assume
without loss of generality that $M$ has rank 3.  Then, to show that
$B(\tau)$ is flat, it suffices to find a faithful representation which
deforms flatly.  The local system $\pi_!(\CC)$ affords such a
representation, since by assumption $\CC_\psi[W_+]$ deforms flatly to
$\hat{A}_+(M)/J\hat{A}_+(M)$.
\end{proof}

\subsection{End of proof of Theorem \ref{formal pbw}.}
Now we can finish the proof of Theorem \ref{formal pbw}.  We can regard
${\mathbf M}$ as a $B$-module (in which all the arrows are isomorphisms).  By
Lemma \ref{cohom}, ${\rm Ext}^1_B({\mathbf M},{\mathbf M})={\rm Ext}^2_B({\mathbf
  M},{\mathbf M})=0$.  By standard deformation theory, this implies that
${\mathbf M}$ can be uniquely deformed to a module ${\mathbf M}(\tau)$ over
$B_J(\tau)$.  This module still has the property that all maps are
isomorphisms.  Thus the $\hat{A}_{ijk+}/J\hat{A}_{ijk+}$-module
structure on the subspace ${\mathbf M}(\tau)_{ijk}$ transports to such a
structure on the subspace ${\mathbf M}(\tau)_N$ (the fiber at the
north pole), and all of the different
module structures will be compatible.  In other words, ${\mathbf M}(\tau)_N$
is a representation of $\hat{A}_+/J\hat{A}_+$ that flatly deforms the
regular representation of $\CC_\psi[W_+]$. The existence of such a
deformation implies the flatness of $\hat A_+/J\hat A_+$. Theorem
\ref{formal pbw} is proved.

\section{A generalization of the Iwahori-Hecke algebra}

Although Theorem \ref{mai} is stated in terms of deformations of even
subgroups of Coxeter groups, it also gives rise to interesting deformations
of Coxeter groups themselves.  Indeed, if $u\in \mathbb T_f$ is invariant
under the $\ZZ_2$ action (modulo the permutation action), then
$\tilde{A}_{\pi(u)}$ is an algebra deforming $\ZZ[W]$.  The corresponding
flatness locus is obtained from $\mathbb T_f$ by adding the additional
``edge'' conditions
\[
(e^{(m_{ij})}_{ij})^2 = 1,\ e^{(k)}_{ij} = e^{(m_{ij})}_{ij}e^{(m_{ij}-k)}_{ij}.
\]

Of particular interest is the case $e^{(m_{ij})}_{ij}=(-1)^{m_{ij}-1}$,
$e^{(m_{ij}/2)}_{ij}=0$ (the latter being redundant when $2$ is
invertible).  It turns out that in this case, the flatness conditions are
greatly simplified: only the $(2,3,3)$ case gives an additional condition,
namely that if $m_{ij}=m_{jk}=3$, $m_{ik}=2$, then
\[
e^{(1)}_{ij}=e^{(1)}_{jk}.
\]
This gives rise to the following generalization of the Iwahori-Hecke
algebra.  If $x$, $y$ are elements of an associative algebra, let
$B_k(x,y)$ denote the corresponding braid relation; that is,
\[
B_{2k}(x,y) = (xy)^k-(yx)^k,\quad
B_{2k+1}(x,y) = (xy)^kx-(yx)^ky.
\]

\begin{thm}
Introduce $u_i$, $v_i$, $1\le i\le r$, and for $i<j$ such that
$m_{ij}<\infty$, introduce parameters $f^{(l)}_{ij}$, $1\le l<m_{ij}/2$,
such that the following two conditions hold:
\begin{itemize}
\item[(1)] If $m_{ij}<\infty$ is odd, then $u_i=u_j$, $v_i=v_j$.
\item[(2)] If $m_{ij}=m_{jk}=3$, $m_{ik}=2$, then
  $f^{(1)}_{ij}=f^{(1)}_{jk}$.
\end{itemize}
Let $A'$ be the $\ZZ[u_i,v_i,f^{(l)}_{ij}]$-algebra with presentation
\begin{align}
\langle
T_1,\dots,T_r
|&
T_i^2-u_i T_i+v_i,\notag\\
&B_{m_{ij}}(T_i,T_j)
+
f^{(1)}_{ij} B_{m_{ij}-2}(T_i,T_j)
+
f^{(2)}_{ij} B_{m_{ij}-4}(T_i,T_j)
+\cdots
\rangle,\notag
\end{align}
Then $A'$ is a free $\ZZ[u_i,v_i,f^{(l)}_{ij}]$-module, with basis of the
form $T_w$, $w\in W$.  In particular, $A'$ is a flat deformation of $\ZZ[W]$.
\end{thm}

\begin{proof}
If we choose a reduced word for each element $w\in W$, the corresponding
products $T_w$ certainly span $A'$, and it suffices to show that they are
linearly independent.  To show this,  we may adjoin
$1/2$ and $1/\sqrt{u_i^2-4v_i}$ for all $i$ to the coefficient ring.  At
this point, if we consider new generators
$U_i=(2T_i-u_i)/\sqrt{u_i^2-4v_i}$ we need simply show that the
corresponding $U_w$ are linearly independent (since the two sets are
related by a triangular change of basis with unit diagonal).  In terms of
the $U_i$, the deformed braid relations have the same form (given condition
(1) above), and thus it suffices to prove linear independence when $u_i=0$,
$v_i=1$.

At this point, if we multiply the deformed braid relation by
$(T_jT_i)^{m_{ij}/2}$ or $T_j(T_iT_j)^{(m_{ij}-1)/2}$ as appropriate, we
obtain a polynomial in $T_iT_j$.  But the resulting algebra is then a
specialization $\tilde{A}_{\pi(u)}$ satisfying the invariance and flatness
conditions.
\end{proof}

In particular, the case $f^{(l)}_{ij}\equiv 0$ is the ordinary
Iwahori-Hecke algebra associated to $W$.

\begin{remark}
Other 2-cocycles of $W$ will of course have other deformations associated to
them; the difficulty, however, is that the translation by $u_i/2$ no longer
preserves the form of the deformed braid relation in those cases.
\end{remark}

Another way of obtaining flat deformations of $\ZZ[W]$ is via the
observation that every Coxeter group is also the even subgroup of a
slightly larger Coxeter group: $W = (W\times \ZZ_2)_+$.  The resulting
algebras include all algebras of the form $\tilde{A}_{\pi(u)}$, as well as
deformations in which the quadratic relations are more general.  It is
unclear, however, whether these deformations are truly more general; it is
possible that the extra freedom can be removed as in the above Hecke
algebra case.

\section{Additive versions}

In this section, for simplicity, we work over $\mathbb C$ rather
than over $\mathbb Z$. 

\subsection{Graded algebras} 

Recall \cite{ER} that the polynomial relation for the algebra $A$ 
can be written as 
$$
s_is_j...=(-1)^{m_{ij}+1}t_{ij}s_js_i...+S.L.T.
$$
($m_{ij}$ factors on both sides), where $t_{ij}:=
\prod_k t_{ijk}$, and S.L.T. stands for ``smaller length terms''. 
Thus it is natural to define the graded version $\oA$ of the algebra $A$, 
generated over $\overline{R}:=\mathbb C[t_{ij}]$ ($t_{ji}=t_{ij}^{-1}$) 
by $s_i$ with defining relations 
$$
s_i^2=0,\ s_is_j...=(-1)^{m_{ij}+1}t_{ij}s_js_i...,\quad s_it_{jp}=t_{jp}^{-1}s_i.
$$
Let $\oA_t$ be the fiber of $\oA$ over $t\in \overline{\mathbb T}:={\rm
Spec}\overline{R}$ (not always an algebra). 

It is easy to see that $T_{w(x)}$, $x\in W$, is still a 
spanning set of $\oA$ over $\overline{R}$ and of $\oA_t$ over 
$\mathbb C$, and that if $t\in \Theta$ then $T_{w(x)}$ are a basis
of $\oA_t$. Also, if 
$T_{w(x)}$ are a basis, the Hilbert series
of $\oA_t$ is equal to $h(z)=h_M(z)$.

Let us now consider the algebra $\oA_1$, corresponding to the
special case $t=1$. Note that by Lemma \ref{spin}, $1\in \Theta$, so 
in this case $T_{w(x)}$ are a basis. 

\begin{lemma}\label{addmult}
(i) The algebra $A_1$ is generated by $s_i$ with defining
relations $s_i^2=1$, $(s_i-s_j)^{m_{ij}}=0$.

(ii) The algebra $\oA_1$ is generated by $s_i$ with defining
relations $s_i^2=0$, $(s_i-s_j)^{m_{ij}}=0$.
\end{lemma}

\begin{proof}
(i) We will show that the relation $(s_is_j-1)^{m_{ij}}=0$ is equivalent to
  $(s_i-s_j)^{m_{ij}}=0$ in the presence of the relations
  $s_i^2=s_j^2=1$. To do so, we write
$$
(s_i-s_j)^{m_{ij}}=(s_is_j-1)s_j*s_i(1-s_is_j)*(s_is_j-1)s_j...=\pm
(s_is_j-1)^{m_{ij}}s_js_is_j...,
$$
since $s_js_i$ is the inverse of $s_is_j$ and hence commutes with
it. This proves (i).

(ii) The result follows by opening brackets in the relation 
$(s_i-s_j)^{m_{ij}}=0$, and writing down the resulting
$2^{m_{ij}}$ terms.
\end{proof}

\subsection{Additive versions} 

Let us define the ``additive'' versions of 
the algebras $A,A_+$.
Namely, let $\tau_{ijk}$ be parameters, such that
$\tau_{ijk}=-\tau_{ji,-k}$, $k\in \mathbb Z_{m_{ij}}$
($m_{ij}<\infty$). 

Define the algebra ${\mathcal A}_+$ generated over ${\mathcal
R}:=\mathbb C[\tau]$ 
by generators $\alpha_{ij}$, $i\ne j$, with defining relations 
$$
\alpha_{ij}+\alpha_{ji}=0,\ \alpha_{ij}+\alpha_{jk}+\alpha_{ki}=0,\ 
\prod_{k=1}^{m_{ij}}(\alpha_{ij}-\tau_{ijk})=0
\ (m_{ij}<\infty).
$$ 
Further, define ${\mathcal A}$ to be the semidirect product $\mathbb C\mathbb
Z_2\ltimes {\mathcal A}_+$, where the nontrivial element 
$\sigma$ of $\mathbb Z_2$ acts as follows: 
$\sigma(\alpha_{ij})=-\alpha_{ij}$, $\sigma(\tau_{ijk})=\tau_{jik}$. 
We'll denote by ${\mathcal A}_0$, ${\mathcal A}_{0+}$ 
the zero-fibers of these algebras. Thus, ${\mathcal A}_{0+}$
is generated by $\alpha_{ij}$ with defining relations 
$$
\alpha_{ij}+\alpha_{ji}=0,\ \alpha_{ij}+\alpha_{jk}+
\alpha_{ki}=0,\ \alpha_{ij}^{m_{ij}}=0\
(m_{ij}<\infty).
$$ 
This algebra is  naturally graded by $\mathbb Z_+$ ($\deg(\alpha_{ij})=1$).

\begin{proposition}\label{spset}
Let $0\in I$. For any element $x\in W$ such 
that $l(s_0x)>l(x)$, fix a reduced decomposition 
$x=s_{i_1}...s_{i_n}$. Set $b_x:=
\alpha_{0i_1}\alpha_{i_1i_2}...\alpha_{i_{n-1}i_n}$
(we agree that $b_1=1$). 
Then the collection of all the elements $b_x$ 
is a spanning set of ${\mathcal A}_{0+}$. In particular, the 
Hilbert series of ${\mathcal A}_{0+}$ is dominated
(coefficientwise) by $h(z)/(1+z)$. 
\end{proposition}

\begin{proof} Since $\alpha_{ij}+\alpha_{ji}=0,
\alpha_{ij}+\alpha_{jk}+
\alpha_{ki}=0$, we see that 
any element of ${\mathcal A}_{0+}$ is a linear combination 
of elements of the form 
$\alpha_{0i_1}\alpha_{i_1i_2}...\alpha_{i_{n-1}i_n}$
for some $i_1,...,i_n$. 

Now let us apply the following ``braid'' relations of degree
$m_{ij}$, which are clearly satisfied in ${\mathcal A}_{0+}$:
$$
\alpha_{ij}\alpha_{ji}...\alpha_{ji}\alpha_{ip}=
\alpha_{ji}\alpha_{ij}...\alpha_{ij}\alpha_{jp},
$$
$$
\alpha_{pi}\alpha_{ij}\alpha_{ji}...\alpha_{ji}=
\alpha_{pj}\alpha_{ji}\alpha_{ij}...\alpha_{ij}
$$
if $m_{ij}$ is odd, and 
$$
\alpha_{ij}\alpha_{ji}...\alpha_{ij}\alpha_{jp}=
-\alpha_{ji}\alpha_{ij}...\alpha_{ji}\alpha_{ip},
$$
$$
\alpha_{pi}\alpha_{ij}\alpha_{ji}...\alpha_{ij}=
-\alpha_{pj}\alpha_{ji}\alpha_{ij}...\alpha_{ji}
$$
if $m_{ij}$ is even. Using these braid relations, we can 
reduce the monomial $\alpha_{0i_1}\alpha_{i_1i_2}...\alpha_{i_{n-1}i_n}$
to zero if the word $s_0....s_n$ is not reduced, 
and to the monomial $b_x$ for some $x$ if this word is reduced. 
The proposition is proved. 
\end{proof}

\begin{proposition}\label{homom}
(i) There exists a (unique) homomorphism $\phi_0: {\mathcal A}_{0+}\to
\oA_1$ defined by the formula $\phi_0(a_{ij})=s_i-s_j$.

(ii) There exists a (unique) homomorphism $\phi: {\mathcal A}_{0+}\to
A_1$ defined by the formula $\phi(a_{ij})=s_i-s_j$.
\end{proposition}

\begin{proof}
This follows immediately from Lemma \ref{addmult}.
\end{proof}

\begin{lemma}\label{sum}
Let $B=Im\phi_0$, i.e. the subalgebra of $\oA_1$
generated by $s_i-s_j$. Then for any $0\in I$, $\oA_1=B+s_0B$. 
\end{lemma}

\begin{proof} 
We first show that any element of $\oA_1$ can be written 
as a linear combination $x+s_0y$, $x,y\in B$.
It suffices to consider elements $z=s_{i_1}...s_{i_n}$, 
and use induction in $n$. If $n=0,1$, the statement is clear
($s_i=s_0+(s_i-s_0)$). If $n>1$, we have 
$$
z=s_{i_1}...s_{i_{n-1}}s_{i_n}=s_{i_1}s_{i_2}...s_{i_{n-1}}(s_{i_n}-
s_{i_{n-1}}).
$$ 
Since $s_1...s_{n-1}\in B+s_0B$ by the induction assumption, we
get $s_1...s_n\in B+s_0B$. 

Thus, $\oA_1=B+s_0B$.
\end{proof} 

The main result of this section is the following theorem.

\begin{theorem}\label{bas}
(i) The map $\phi_0: {\mathcal A}_{0+}\to B$ is an isomorphism.

(ii) The Hilbert series of ${\mathcal A}_{0+}$ and $B$
is $h(z)/(1+z)$. 

(iii) $\oA_1=B\oplus s_0B$. 

(iv) The spanning set of ${\mathcal A}_{0+}$ from Proposition
\ref{spset} is in fact a basis.  
\end{theorem}

\begin{proof} By Lemma \ref{sum}, the Hilbert series 
of $B$ is at least $h(z)/(1+z)$, therefore so is the Hilbert
series of ${\mathcal A}_{0+}$. On the other hand, by Proposition 
\ref{spset}, the opposite inequality holds. This implies that the
Hilbert series in question is exactly $h(z)/(1+z)$. (i)-(iv)
immediately follow. 
\end{proof}

We conjecture that (ii) remains true for the fibers 
${\mathcal A}_{\tau+}$ over the locus of
$\tau$ such that (ii) holds for finite triangles; equivalently
for $\tau$ in the tangent cone to ${\mathbb T}_f$ at $1$.
This conjecture is true in the rank 3 case by using the Riemann-Hilbert
correspondence between solutions of the additive and
multiplicative Deligne-Simpson problems.

\end{document}